\documentclass[draft]{amsart}
\usepackage{amsmath,amsthm,amssymb}

\newtheorem{theorem}{Theorem}
\newtheorem{prop}{Proposition}
\newtheorem{lem}{Lemma}
\newcommand{\C}{\mathbb{C}\mkern1mu}

\newcommand{\mf}{\mathfrak}
\begin{document}
\title{Antiholomorphic involutions of spherical complex spaces}
\author{Dmitri Akhiezer}
\address{Institute for Information Transmission Problems, B.Karetny 19, 101447 Moscow, Russia}
\email{akhiezer@mccme.ru}
\author{ Annett P\"uttmann }
\address{Ruhr-Universit\"at Bochum, Fakult\"at f\"ur Mathematik,
Universit\"atsstra\ss e 150, 44780 Bochum, Germany}
\email{annett.puettmann@rub.de}
\date{}
\thanks{Research supported by SFB/TR12 ``Symmetrien und Universalit\"at in mesoskopischen 
Systemen'' of the Deutsche Forschungsgemeinschaft. First author supported in part
by the Russian Foundation for Basic Research, Grant 04-01-00647}
\begin{abstract}
Let $X$ be a holomorphically separable irreducible reduced complex space,
$K$ a connected compact Lie group acting on $X$ by holomorphic transformations,
$\theta : K \to K$ a Weyl involution, and $\mu : X \to X$
an antiholomorphic map satisfying $\mu ^2 ={\rm Id}$ and $\mu (kx) = \theta (k)\mu (x)$
for $x\in X,\ k\in K$. We show that if ${\mathcal O}(X)$ is a multiplicity free $K$-module
then $\mu $ maps every $K$-orbit onto itself.
For a spherical affine homogeneous space $X=G/H$ of
the reductive group $G=K^{\mathbb C}$ we construct an antiholomorphic map $\mu $  with these properties.
\end{abstract}
\maketitle
\section{Introduction}
Let $X=(X,\mathcal{O})$ be a complex space on which a compact Lie group $K$ acts continuously by holomorphic transformations. 
Then the Fr\'echet space $\mathcal{O}(X)$ has a natural structure of a $K$-module.
Recall that a $K$-module $W$ is called {\it multiplicity free} if any irreducible $K$-module occurs in $W$ 
with multiplicity 1 or does not occur at all. 
A self-map $\mu $ of a complex space $X$ is called an {\it antiholomorphic involution} 
if $\mu $ is antiholomorphic and $\mu ^2 = {\rm Id}$.
For complex manifolds, J.Faraut and E.G.F.Thomas gave an interesting and simple geometric condition
which implies that $\mathcal{O}(X)$ is a multiplicity free
$K$-module, see \cite {FT}. Namely, for a complex manifold $X$ the $K$-action in $\mathcal{O}(X)$ is multiplicity free if

\medskip
\noindent 
(FT) there exists an antiholomorphic involution $\mu : X \to X$ with the property that, 
 for every $x\in X$, there is an element $k\in K$ such that $\mu (x) = k\cdot x$.

\medskip
\noindent
The proof of Theorem 3 in \cite {FT} goes without changes for irreducible reduced
complex spaces. It should be noted that the setting in \cite {FT} is more general. Namely, the authors consider
any, not necessarily compact, group of holomorphic transformations of $X$
and study invariant Hilbert subspaces of ${\mathcal{O}}(X)$. We will give a simplified proof of their result in our context,
see Proposition \ref{FTshort}. Our main purpose, however, is to prove the converse theorem for a special class of manifolds, namely, for
Stein (or, equivalently, affine algebraic) homogeneous spaces of complex reductive groups. 

Let $G$ be a connected reductive complex algebraic group and $K\subset G$ a maximal compact subgroup.
We prove
that an affine homogeneous space $X = G/H$ is spherical (or, equivalently, 
${\mathcal O}(X)$ is 
a multiplicity free $K$-module) if and only if the $K$-action on $X$ satisfies (FT),
see Theorem \ref{thm:affine}.    
Recall that a diffeomorphism $\mu $ of a manifold $X$ with a $K$-action is said to be {\it $\theta $-equivariant},
if $\theta $ is an automorphism of $K$ and 
$\mu (kx) = \theta (k)\mu (x)$ for all $x\in X,\ k\in K$.
For $X=G/H$ spherical we can say more about $\mu $ in (FT). Namely,
again by Theorem \ref{thm:affine}, 
$\mu $ can be chosen $\theta $-equivariant, where $\theta $
is a Weyl involution of $K$.

In order to prove Theorem \ref{thm:affine}, we consider $\theta $-equivariant antiholomorphic involutions
in a more general context. Namely, let $X$ be a holomorphically separable
irreducible reduced complex space, $K$ a connected compact Lie group of holomorphic transformations of $X$,
and $\mu $ an antiholomorphic $\theta $-equivariant involution of $X$. 
Then our Theorem \ref{thm:sep} asserts that ${\mathcal O}(X)$ is multiplicity free if and only
if $\mu(x)\in Kx$ for all $x\in X$, i.e., if $X$ has property (FT) with respect to $\mu $.

Another important ingredient in the proof of Theorem \ref{thm:affine}
is the construction of two commuting involutions of $G$, a Weyl involution and a Cartan involution,
which both preserve a given reductive spherical subgroup $H \subset G$, see Theorem \ref{thm:taumu}.
The proof is based on the results of \cite{AV} and, therefore, on the classification
of spherical subgroups.

At the end of our paper we give an example of an affine homogeneous space without $\theta $-equivariant 
antiholomorphic involutions,
see Proposition \ref{inv}.

\section{Fourier series of Harish-Chandra }

Harish-Chandra carried over the classical Fourier series
to the representation theory of compact Lie groups
in Fr\'echet spaces, see \cite{H-C}. In this paper,  we will need only the representations in ${\mathcal O}(X)$,
where $X$ is a complex space. We recall the result of Harish-Chandra in this setting.
The details can be found in \cite{H-C}, see also \cite {A}, Ch. 5. 

Let $K$ be a compact Lie group, $\hat K$ its unitary dual, $dk$ the normalized Haar measure
on $K$. For $\delta \in \hat K$ let $\chi _\delta $ denote
the character of $\delta $ multiplied by the dimension of $\delta $. Suppose that $K$ acts by holomorphic transformations
on a complex space $X$. Then 
we have a continuous representation of $K$ in ${\mathcal C}(X)$ and in ${\mathcal O}(X)$.
We will assume that $X$ is reduced and irreducible, so the representation is given by
$k\cdot f(x) = f(k^{-1}x)$, where $k\in K,x\in X$. 

Define an operator family $\{E_\delta \}_{\delta \in \hat K}$ in ${\mathcal O}(X)$ by
$$E_\delta f (x) = \int _K \ \overline {\chi _\delta (k)}\cdot f(k^{-1}x) \cdot dk.$$ 
From orthogonality relations for characters it follows that
all $E _\delta $ commute with the representation of $K$. Furthermore, $\{E_\delta \}_{\delta \in \hat K}$
is a family of projection operators, i.e.  
$E_\delta ^2 = E_\delta $ and $E_\delta E_\epsilon = 0 $ if $\delta \ne \epsilon $.
Let ${\mathcal O}_\delta (X) = E_\delta {\mathcal O}(X)$. Then ${\mathcal O}_\delta (X) = {\rm Ker}(E_\delta - {\rm Id})$,
so ${\mathcal O}_\delta (X)$ is a closed subspace. Again from orthogonality relations it follows that
${\mathcal O}_\delta (X)$ is the {\it isotypic component of type $\delta $}, i.e., ${\mathcal O}_\delta (X)$
consists of all those vectors in ${\mathcal O}(X)$ whose $K$-orbit is contained in
a finite-dimensional $K$-submodule where the representation is some multiple of $\delta $.

Harish-Chandra proved that
each $f\in\mathcal{O}(X)$ can be uniquely written in the form  
$$f=\sum_{\delta \in \hat K}\ f_{\delta},$$ 
where 
$f_{\delta} = E_\delta (f) \in\mathcal{O}_{\delta}(X)$ 
and the convergence is absolute and uniform on compact subsets in $X$.

Assume now that $L$ is another compact Lie group acting by holomorphic transformations
of another complex space $Y$ subject to our assumptions. We will use similar notation for $L$,
in particular, $\theta _\epsilon $ will denote the character of $\epsilon \in \hat L$
multiplied by the dimension of $\epsilon $.
For the representation of $K\times L$ in ${\mathcal O}(X\times Y)$ defined by
$$(k,l) \cdot f(x,y) = f(k^{-1}x, l^{-1}y),\ x\in X, y\in Y, k\in K, l\in L,$$
the type of an isotypic component is determined by a pair $\delta \in \hat K,
\epsilon \in \hat L$. The corresponding isotypic component will be denoted
${\mathcal O} _{\delta,\epsilon}(X\times Y)$.
Of course, the tensor product ${\mathcal O}_{\delta}(X) \otimes {\mathcal O}_{\epsilon}(Y)$
is contained in 
${\mathcal O}_{\delta,\epsilon}(X\times Y)$.
We will need the following lemma.
\begin{lem}
\label{tensor}
If ${\rm dim}\ {\mathcal O}_{\delta }(X) < \infty $ for some $\delta \in \hat K$ then
${\mathcal O}_{\delta,\epsilon}(X\times Y) = {\mathcal O}_{\delta}(X) \otimes {\mathcal O}_{\epsilon}(Y)$
for all $\epsilon \in \hat L$.
\end{lem}
\begin{proof}
Let $f\in {\mathcal O}_{\delta, \epsilon}(X\times Y)$, then
\begin{align*}
f(x,y) & = \int _{K\times L}\ \overline {\chi _\delta (k) \theta _\epsilon (l)}
\cdot f(k^{-1}x, l^{-1}y) \cdot dk dl = \\
& =\int _L \overline {\theta _\epsilon (l )} \, \Bigl( \int _K\ 
\overline{\chi _\delta (k)} \cdot f(k^{-1}x, l^{-1}y)\cdot dk \Bigr) \, dl
\end{align*}
by Fubini theorem.
The function
$$ x\mapsto \int _K\ \overline {\chi _\delta (k)}\cdot f(k^{-1}x,y)\cdot dk $$ 
is in ${\mathcal O}_{\delta} (X)$ for all $y \in Y$. Let $\{\varphi _i\} _{i=1,\ldots,N }$ be a basis
of ${\mathcal O}_\delta (X)$. Then
$$ \int _K\ \overline {\chi _\delta (k)}\cdot f(k^{-1}x,y)\cdot dk = \sum  _{i=1}^N  \  c_i(y)\varphi _i(x)$$
with some $c_i\in {\mathcal O}(Y)$.
Replace in this equality $y$ by $l^{-1}y$, multiply it by $\overline {\theta _{\epsilon}(l)}$
and integrate over $L$ against the Haar measure $dl$. Then we get
$$f(x,y) = \sum _{i=1}^N\ \varphi _i(x)\psi _i(y),$$
where
 $$\psi _i(y) = \int _L \ \overline {\theta _\epsilon (l)}\cdot c_i(l^{-1}y)\cdot dl \ \in {\mathcal O}_\epsilon (Y).$$ 
\end{proof}

\section{$K$-action and complex conjugation}

As in the previous section,
$X$ is an irreducible reduced complex space and $K$ is a compact group acting on
$X$ by holomorphic transformations.
\begin{lem}
\label{fin}
Let  
$W\subset\mathcal{O}(X)$ be a finite-dimensional $K$-submodule.
Introduce a $K$-invariant Hermitian inner product and choose a unitary basis
$\{ f_1,\ldots,f_N\}$ in $W$. 
The function $F:=\sum_{j=1}^N f_j\overline{f_j}$ is $K$-invariant.
Furthermore, $F$ does not depend on the choice of basis.
\end{lem}
\begin{proof}
Let $\{g_1,\ldots,g_N\}$ be another unitary basis of $W$.
There is a unitary transformation $A:W\to W $ such that 
$A(f_j) = g_j = \sum_{l=1}^N a_{lj}f_l$.
We have 
$$\sum_{j=1}^N g_j\overline{g_j} = 
\sum_{j=1}^N \sum_{l,l'=1}^N a_{lj}\overline{a_{l'j}}f_l\overline{f_{l'}}
= \sum_{l,l'=1}^N \delta_{ll'} f_l\overline{f_{l'}}=F.$$
Now, $k\cdot F = \sum_{j=1}^N (k\cdot f_j)\overline{(k\cdot f_j)}$ for $k\in K$.
But $\{ k\cdot f_1,\ldots,k\cdot f_N\}$ is another unitary basis of $W$.
Since $F$ does not depend on the choice of basis, 
it follows that $k\cdot F=F$ for any $k\in K$.
\end{proof}
\begin{lem}
\label{dual}
If $W\subset\mathcal{O}(X)$ is a finite-dimensional $K$-submodule,
then $\overline{W}\subset \overline{\mathcal{O}(X)}$ is also a $K$-submodule,
which is isomorphic to the dual module $W^*$.  
\end{lem}
\begin{proof}
Let $(f,g)$ be a
$K$-invariant Hermitian product on $W$.
For $f\in W,\ \phi \in \overline{W}$ we have the bilinear pairing
$$\langle f,\phi \rangle
= (f,\overline{\phi}),$$
which is obviously $K$-invariant and non-degenerate. This shows that $\overline {W}$
is isomorphic to $W^*$.
\end{proof}
\noindent
Let 
$\mu : X \to X$ be an antiholomorphic involution. Then, by definition, the function $\mu f (x) = f(\mu x)$
is antiholomorphic for any $f \in {\mathcal O}(X)$. 
We want to give a simple proof of the theorem of J.Faraut and E.G.F.Thomas in our setting.
\begin{prop}
\label{FTshort}
If the $K$-action on $X$ satisfies (FT) then ${\mathcal O}(X)$
is a multiplicity free $K$-module.
\end {prop}
\begin{proof} 
Assume the contrary. Let $W, W^\prime \subset {\mathcal O}(X)$ be two irreducible isomorphic $K$-submodules, such that
$W \ne W^\prime $.
Define $f_1,\ldots,f_N\in W$ as in Lemma \ref{fin}.
Fix a $K$-equivariant isomorphism $\phi : W \to W^\prime $ and let $f_i^\prime = \phi (f_i)$.
By Lemma \ref{fin} the function $F = \sum f_i\overline {f_i}$ is $K$-invariant.
The same proof shows that the function $G = \sum f_i\overline{f_i^\prime}$ is also $K$-invariant.
By (FT) we have $\mu F = F$ and $\mu G = G$.
Since the multiplication map ${\mathcal O}(X) \otimes \overline {{\mathcal O}(X)}\to {\mathcal O}(X)\cdot\overline{{\mathcal O}(X)}$
is an isomorphism of vector spaces,
it follows that
$$\sum _i \overline {\mu f_i}\otimes \mu f_i = \sum _i f_i\otimes \overline{f_i}$$
and
$$\sum _i \overline{\mu f_i^\prime}\otimes \mu f_i = \sum _i f_i\otimes \overline {f_i^\prime}.$$
Therefore
the linear span of $\overline {\mu f_1},\ldots,\overline {\mu f_N}$ coincides with the linear span of
$f_1, \ldots, f_N$ 
and with the linear 
span of $\overline{\mu f_1^\prime},\ldots,\overline{\mu f_N^\prime}$.
 Thus $\overline {\mu W}= \overline {\mu W^\prime }$ and $W = W^\prime $, contradictory to our assumption.  
\end{proof}
From now on the compact group $K$ is assumed connected. An involutive automorphism
$\theta:K\to K$ is called a 
{\it Weyl involution} if there exists a maximal torus $T\subset K$
 such that $\theta(t)=t^{-1}$ for $t\in T$. It is known that Weyl involutions exist and that they are all
conjugated by inner automorphisms. If $\theta $ is a Weyl involution and $\rho $ a linear representation of $K$
then $\rho \circ \theta$ is the dual representation.
Recall that an antiholomorphic involution $\mu : X \to X$ is called $\theta$-equivariant
if $\mu (kx) = \theta (k)\mu (x)$ for all $x\in X,\ k\in K$.
 
\begin{lem}
\label{anti}
Let $\theta $ be a Weyl involution of $K$ and $\mu:X\to X$ a $\theta $-equivariant antiholomorphic involution of $X$.
If $W\subset\mathcal{O}(X)$ is a finite-dimensional $K$-submodule then
$\mu W$ is also a $K$-submodule. Furthermore, $\overline{W}$ and $\mu W$ are isomorphic $K$-modules.
\end{lem}
\begin{proof}Introduce a $K$-invariant Hermitian inner product and choose a unitary basis
$\{ f_1,\ldots,f_N\}$ in $W$. 
Denote the representation in $W$ by $\rho$. 
The condition $\mu(kx)=\theta(k)\mu(x)$ implies that
$$k\cdot \mu f(x)=\mu f(k^{-1}x)=f(\mu(k^{-1}x))=f(\theta(k)^{-1}\mu(x))=\theta(k)f(\mu x) = \mu \theta (k)f (x). $$
Hence $\mu W$ is indeed a $K$-submodule with the representation $\rho\circ\theta$.
Since $\theta $ is a Weyl involution, this representation is dual to $\rho $.
But the representation in $\overline {W}$ is also dual to $\rho $ by Lemma \ref{dual},
and our assertion follows.
\end{proof}

\begin{lem}
\label{anti1}
Keep the notation of Lemma \ref {anti} and assume in addition that $W$ is irreducible and $\mu W = \overline {W}$.
Then for a $K$-invariant Hermitian inner product on $W$ one has
 $$(\overline {\mu f_1}, \overline {\mu f_2}) = (f_1,f_2),$$
where $f_1,f_2\in W$.

\end{lem}
\begin{proof}
The new Hermitian inner product $ \{f_1,f_2\} := (\overline {\mu f_1}, \overline {\mu f_2})$
on $W$ is also $K$-invariant. Since $W$ is an irreducible $K$-module, it follows that
$\{f_1,f_2\} = c (f_1,f_2)$, where $c > 0$. But then
$$ \{\overline {\mu f_1}, \overline {\mu f_2}\} = c (\overline {\mu f_1}, \overline {\mu f_2})$$
and, on the other hand,
$$ \{\overline {\mu f_1}, \overline {\mu f_2}\} = (f_1,f_2)$$
because $\mu $ is an involution. Thus
$$c (\overline {\mu f_1}, \overline {\mu f_2}) = (f_1,f_2) = c^{-1}(\overline {\mu f_1}, \overline 
{\mu f_2}),$$
hence $c^2 = 1$ and $c = 1$.
\end{proof}

\section{ Holomorphically separable spaces}

Since we assume that $K$ is connected, the irreducible representations of $K$
are
determined by their highest weights. 
We denote by $W_\lambda $ an irreducible $K$-module
with highest weight $\lambda $ and write ${\mathcal O}_\lambda (X)$ instead of ${\mathcal O}_\delta (X)$,
where $\delta \in \hat K$ and $\lambda = \lambda (\delta )$ is the highest weight of $\delta $. 
Those highest weights $\lambda $, for which $W_\lambda $ occurs in our $K$-module ${\mathcal O}(X)$,
form an additive semigroup, to be denoted by $\Lambda (X)$. In other words,
$\Lambda (X) $ is the set of highest weights such that ${\mathcal O}_\lambda (X) \ne \{0\}$.   
The subspace of fixed vectors
of a $K$-module $W$ is denoted by $W^K$. 
We remark that if $A$ is an algebra on which $K$ acts as a group of automorphisms then $A^K$
is a subalgebra of $A$. 

\begin{theorem}
\label{thm:sep}
Let $X$ be a holomorphically separable irreducible reduced complex space, $K$ a connected compact Lie group acting on 
$X$ by holomorphic transformations, $\theta:K\to K$ a Weyl involution
and $\mu:X\to X$ a $\theta $-equivariant antiholomorphic involution of $X$.
Then $\mathcal{O}(X)$ is multiplicity free if and only if $\mu(x)\in Kx$ for all $x\in X$.
\end{theorem}
\begin{proof}
If $\mu(x)\in Kx$ for all $x\in X$, then (FT) guarantees that
the $K$-action on $\mathcal{O}(X)$ is multiplicity free, 
see Introduction and Proposition \ref{FTshort}.

We now prove the converse.
Let $\mathcal{A}(X) = \mathcal{O}(X)\cdot \overline{\mathcal{O}(X)}$.
Since $X$ is holomorphically separable, the algebra $\mathcal{A}(X)$ separates points of $X$.
By Stone-Weierstrass theorem $\mathcal{A}(X)$ is dense in the algebra $\mathcal{C}(X)$
of continuous functions on $X$. The standard averaging argument shows that
$\mathcal{A}(X)^K$ is dense in $\mathcal{C}(X)^K$. Now, 
if $Kx$ and $Ky$ are two different $K$-orbits in $X$ then there is a $K$-invariant continuous 
function $f\in\mathcal{C}(X)$ which seperates these orbits.
Since this function can be approximated
by $K$-invariant functions from $\mathcal{A}(X)$, it follows that $\mathcal{A}(K)^K$
separates $K$-orbits.
Let $\lambda \in \Lambda (X)$ be a highest weight which occurs in the decomposition of 
the $K$-algebra $\mathcal{O}(X)$. Since $\mathcal{O}(X)$ is multiplicity free, the isotypic
component $\mathcal{O}_\lambda (X)$ is irreducible. We can identify 
this isotypic component with $W_\lambda $, and so we write
$W_\lambda = \mathcal{O}_\lambda(X)$. 
Now apply Lemma \ref{fin} to construct a $K$-invariant function in $W_\lambda
\cdot \overline{W_\lambda}$. Call this function $F_\lambda $.
We claim that the family $\{F_\lambda\} _{\lambda \in \Lambda (X)}$
also separates $K$-orbits in $X$. 

To prove the claim it is enough to present each $F\in \mathcal{A}(X)^K$
as the sum of a series
$$F = \sum _{\lambda \in
  \Lambda (X)}\ c_\lambda F_\lambda\ , $$
where the convergence is absolute and uniform on compact subsets in $X$.
In order to prove this decomposition, consider the complex space $\overline{X}$ with the 
conjugate complex structure. 
There is a natural
$K$-action on $\overline {X}$, and so we obtain an action of $K\times K$
on $X\times \overline{X}$. Since ${\mathcal O}(\overline {X}) = \overline {{\mathcal O}(X)}$,
the isotypic components of the $K$-module ${\mathcal O}(\overline {X})$ are just the submodules
$\overline {W_\lambda}$. 
By Lemma \ref{tensor} the isotypic components of the $(K\times K)$-module ${\mathcal O}(X\times \overline {X})$
are the tensor products $W_\lambda \otimes \overline {W_{\lambda ^\prime}}$.

For any $F\in {\mathcal O}(X\times \overline {X})$
the theorem of Harish-Chandra yields the decomposition
$$F= \sum F_{{\lambda \lambda ^\prime}}  \ \ {\rm with} \ F_{\lambda \lambda ^\prime } \in W_\lambda 
\otimes \overline {W_{\lambda ^\prime}},$$
where the convergence is absolute and uniform on compact subsets in $X\times \overline {X}$.
In particular,
if $F\in(\mathcal{O}(X)\otimes \overline{\mathcal{O}(X)})^K$ 
then all summands are $K$-invariant. But 
$\overline {W_\lambda}$ is dual to $W_\lambda$ by Lemma \ref{dual}, hence $F_{\lambda \lambda ^\prime }=0$ for $\lambda ^\prime
 \ne \lambda $ by Schur lemma.
The remaining summands $F_{\lambda \lambda }$ are 
$K$-invariant elements in $W_\lambda \otimes \overline {W_\lambda }$. But
the space $(W_\lambda \otimes \overline {W_\lambda })^K$ is one-dimensional,
again by Schur lemma. Therefore, restricting $F_{\lambda \lambda}$ 
to the diagonal in $X\times \overline {X}$,
we get the functions proportional to the $F_ \lambda $'s
defined above. 

Now, because ${\mathcal O}(X)$ is multiplicity free, it follows from Lemma \ref{anti} that $\overline{W_\lambda} =  \mu W_\lambda $.
Furthermore, Lemma \ref{anti1} shows that the composition of $\mu $ with complex 
conjugation preserves a $K$-invariant
Hermitian product on $W_\lambda $. Therefore $\mu F_\lambda = F_\lambda $ by Lemma \ref{fin}.
Since the family of functions  
$F_{\lambda}$ separates $K$-orbits, $\mu $ must preserve each of them or, equivalently,
$\mu x\in Kx $ for all $x\in X$.
\end{proof}

\medskip
\noindent
{\it Remark}\, For the torus $T = (S^1)^m$ 
the Weyl involution is given by $\theta (t) = t^{-1}$.
Suppose that $T$ acts on ${\mathbb P}_n$ by $t\cdot (z_0:\ldots : z_n) = 
(\chi_0(t)z_0:\ldots:\chi_n(t)z_n)$ with some characters $\chi _i : T \to S^1,\ i=0,\ldots,n,$ and $\mu : {\mathbb P}_n
\to {\mathbb P}_n $ is given by $\mu (z_0:\ldots:z_n) =
(\overline{z_0}:\ldots: 
\overline{z_n})$. Then $\mu $ is obviously $\theta $-equivariant.
However, if $m < n$ then $\mu $ cannot map each $T$-orbit onto itself.
This shows that holomorphic separability of $X$ in Theorem \ref{thm:sep} is essential.   

\medskip
\noindent
Let $\mf{k} $ be the Lie algebra of $K$ and $\mf{g}= \mf{k}^{\C}$
its complexification.
An irreducible reduced complex space $X$ is called {\it spherical} under the action of a 
compact connected Lie group $K$, if $X$ is normal and there exists a point $x\in X$ such that the tangent space
$T_xX$ is generated by the elements of a Borel subalgebra of 
$\mf{g}$, see\cite{AH}.

\begin{theorem}
\label{thm:Stein}
Let $X$ be a normal Stein space, $K$ a connected compact Lie group acting on $X$
by holomorphic transformations, $\theta : K \to K$ a Weyl involution and $\mu : X \to X$
a $\theta$-equivariant antiholomorphic ivolution of $X$. Then $X$ is spherical
if and only if $\mu (x) \in Kx$ for all $x\in X$.
\end{theorem}
\begin{proof} 
It is known that a normal Stein space $X$ is spherical if and only if
$\mathcal{O}(X)$ is a multiplicity free $K$-module \cite{AH}.
The result follows from Theorem \ref{thm:sep}.
\end{proof}

\section{Weyl involution and Cartan involution}

Throughout this section, except Theorem \ref{thm:affine}, 
the word {\it involution} means an involutive automorphism of a group. This notion 
will be used for
complex algebraic groups and for Lie groups. Let
$G$ be a connected reductive algebraic group over ${\mathbb C}$ and let $K$ be a connected
compact Lie group. 
So far we considered Weyl involutions of $K$,
but they can be also defined for $G$.
Namely,
an involution $\theta:G\to G$ is called a 
{\it Weyl involution} if there exists a maximal algebraic torus $T\subset G$
such that $\theta(t)=t^{-1}$ for $t\in T$.  

\begin{lem}
Let $G$ be a connected reductive complex algebraic group and $K\subset G$ a maximal
compact subgroup. Any Weyl involution $\theta$ of $K$ extends uniquely
to a Weyl involution of $G$.
\end{lem}
\begin{proof}
By Theorem 5.2.11 in \cite{OV} any differentiable automorphism $K\to K$
extends uniquely to a polynomial automorphism $G\to G$.
Let $\theta : K \to K $ be a Weyl involution and
$T\subset K$ a maximal torus of $K$ on which $\theta(t)=t^{-1}$.
Now, the complexification $T^{\C}$ of $T$ is a maximal torus of $G$.
The extension of $\theta$ to $G$, which we again denote by $\theta$, is a
Weyl involution of $G$ because $\theta(t)=t^{-1}$ for all 
$t\in T^{\C}$.
\end{proof}
An algebraic subgroup $H\subset G$ is called {\it spherical} if $G/H$ is a spherical variety,
i.e., if a Borel subgroup of $G$ acts
on $G/H$ with an open orbit.
A reductive algebraic subgroup $H\subset G$ is called {\it adapted} if there exists a Weyl involution
$\theta:G\to G$ such that $\theta(H)=H$ and $\theta|_{H^{\circ}}$ is a Weyl involution of the
connected component $H^{\circ}$.
A similar definition is used for compact subgroups of connected compact
Lie groups.
\begin{prop}
Any spherical reductive subgroup
$H\subset G$ is adapted.
\end{prop}
\begin{proof}
See \cite{AV}, Proposition 5.10.
\end{proof}

\begin{prop}
Let $H\subset G$ be an adapted algebraic subgroup,
$K\subset G$ and $L\subset H$ maximal compact subgroups, and $L\subset K$.
Then $L$ is adapted in $K$.
\end{prop}
\begin{proof}
See \cite{AV}, Proposition 5.14.
\end{proof}
\begin{theorem}
\label{thm:taumu}
Let $G$ be a connected reductive algebraic group and $H\subset G$ 
a reductive spherical subgroup.
Then there exist a Weyl involution $\theta:G\to G$ and a Cartan involution $\tau:G\to G$ such that
$\theta\tau = \tau\theta$, $\theta(H)=H$ and $\tau(H)=H$.
\end{theorem}
\begin{proof}
Let $L$ be a maximal compact subgroup of $H$ and $K$ a maximal compact subgroup of $G$ that contains $L$.
Then $K$ is the fixed point subgroup
$G^{\tau}$ of some Cartan involution $\tau$.
Since $H$ is adapted in $G$, there is a Weyl involution $\theta:K\to K$ such that $\theta(L)=L$.
For any $k\in K$, we have $\theta\tau(k)=\theta(k)=\tau\theta(k)$ by the definition
of $\tau$.
Denote again by $\theta:G\to G$ the unique extension to $G$ of the given Weyl involution of
$K$.
Since $G$ is connected and the relation 
$\theta\tau(g)=\tau\theta(g)$ holds on $K$, it also holds on $G$.
\end{proof}
\begin{theorem}
\label{thm:affine}
Let $X=G/H$ be an affine homogeneous space
of a connected reductive algebraic group $G$.
Let $K$ be a maximal compact subgroup of $G$.
Then $X$ is spherical if and only (FT) is satisfied for the action of $K$ on $X$. 
Moreover, if $X$ is spherical one can choose $\mu $ in (FT) to be $\theta $-equivariant,
where $\theta $ is a
Weyl involution of $K$.
\end{theorem}
\begin{proof}
(FT) implies that ${\mathcal O}(X)$ is multiplicity free
or, equivalently, that $X$ is spherical.
Conversely, assume that $X=G/H$ is a spherical variety.
Since $X$ is affine, $H$ is a reductive subgroup by Matsushima-Onishchik theorem. 
Define $\theta $ and $\tau $ as in Theorem \ref{thm:taumu} 
and put $\mu(g\cdot H)=\theta\tau(g)\cdot H$.
The map $\mu : X \to X$ is well defined because $\theta\tau(H)=H$.
The lift of $\mu $ to $G$ is an antiholomorphic involutive automorphism, so
it is obvious that $\mu $ is an antiholomorphic involution of $X$.
Since $\theta\tau=\tau\theta$, it follows that $\theta(K)=K$.
Therefore, for any $x = gH\in X$ one has
$$\mu(kx)=\theta\tau(kg)
\cdot H = \theta(k)\theta\tau(g)\cdot H = \theta(k)\mu(x)$$
for all $k\in K$. From Theorem \ref{thm:Stein} it follows that $\mu (x) \in Kx$ for all $x\in X$. 
\end{proof}

\section {Non-spherical spaces: an example}

We keep the notation of the previous section.  
For a spherical affine homogeneous space $X=G/H$,
we constructed a
$\theta $-equivariant
antiholomorphic involution $\mu $.
In this section we want to show that the sphericity assumption is essential.

\begin{lem}
\label{last1}
Let $X$ be an irreducible reduced complex space with a holomorphic action of  $G$.
Let $\theta $ be any algebraic automorphism of $G$ preserving a maximal compact subgroup $K\subset G$.
Denote by $\tau $ the Cartan involution with fixed point subgroup $K$. If $\mu $ is an antiholomorphic involution of $X$
satisfying $\mu (kx) = \theta (k) \mu (x)$ for all $x \in X,\ k\in K$, then one has $\mu (gx) =
\theta \tau (g)\mu (x)$ for all $x\in X,\ g\in G$.
\end{lem}
\begin{proof}
For every fixed $x\in X$ consider two antiholomorphic maps $\varphi_x:G\to X$
and $\psi_x:G\to X$, defined by
$\varphi_x(g)=\mu (gx)$ and $\psi_x(g)=\theta \tau (g)\mu (x)$.
Since the required identity holds for $g\in K$,
the maps $\varphi_x$ and $\psi_x$ coincide on $K$. 
But $K$ is a maximal totally real submanifold in $G$, so
$\varphi_x$ and $\psi_x$ must coincide on $G$.
\end{proof}
\begin{lem}
\label{last2}
Let $X = G/H$, where $H\subset G$ is an algebraic reductive subgroup, $\theta $ a Weyl involution of $G$
preserving $K$, and $\mu : X \to X$ an antiholomorphic involution of $X$ satisfying $\mu (kx) = \theta (k)\mu (x)$.
Then $H$ and $\theta (H)$ are conjugate by an inner automorphism of $G$.
\end{lem}
\begin{proof}
Assume first that $\tau (H) = H$. Let $x_0 =e\cdot H$ and $h\in H$. Then $\theta (h)\mu (x_0) = \theta (\tau (\tau(h)))\mu (x_0)=
\mu (\tau (h)x_0) = \mu (x_0)$ by Lemma \ref{last1}. It follows that $\theta (H)$ is the stabilizer of $\mu (x_0)$, so $H$
and $\theta (H)$ are conjugate.

To remove the above assumption, take a maximal compact subgroup $L\subset H$ and a maximal compact subgroup $K_1
\subset G$, such that $L\subset K_1$. Then $K_1 = gKg^{-1}$ for some $g\in G$. 
The fixed point subgroup of the Cartan involution $\tau _1 := {\rm Ad}(g)\tau{\rm Ad}(g)^{-1}$ is exactly $K_1$,
so $\tau _1 $ is the identity on $L$ and, consequently, $\tau _1(H)=H$.    
Let $H_1 :=g^{-1}Hg$, then  
$\tau (H_1) = ({\rm Ad} g)^{-1}\tau _1 ({\rm Ad} g)(H_1) 
=({\rm Ad} g)^{-1} \tau _1(H) =({\rm Ad} g)^{-1}(H) = g^{-1}Hg = H_1$.
Replacing $H$ by $H_1$, we can apply the above argument.
\end{proof}
\begin{prop} Let $G = {\rm SO}_{10}(\mathbb{C})$  and $H\subset G$ the adjoint group of ${\rm SO}_{5}(\mathbb {C})$.
Let $\theta $ be a Weyl involution of the maximal compact subgroup
$K={\rm SO}_{10}(\mathbb{R})$. Then an antiholomorphic
involution of $X=G/H$ cannot be $\theta $-equivariant.
\label{inv}
\end{prop}
\begin{proof} Extend $\theta $ holomorphically to $G$.
In view of Lemma \ref{last2} it suffices to show that $\theta (H)$ and $H$ are not 
conjugate
by an inner automorphism of $G$. Assume that $\theta (H) = g_0Hg_0^{-1}$. Then there is an automorphism $\phi : H \to H$,
such that $\theta (h) = g_0\phi (h)g_0^{-1}$ for $h\in H$. All automorphisms of $H$ are inner, so $\phi (h) = h_0hh_0^{-1}$
for some $h_0\in H$. Therefore $\theta (h) = g_1hg_1^{-1}$ for all $h\in H$, where $g_1 = g_0h_0$.
Define an automorphism of $\alpha : G \to G$ by $\alpha : = ({\rm Ad}(g_1))^{-1}\cdot \theta $
and remark that $H \subset G^\alpha $, where $G^\alpha $ is the fixed point subgroup of $\alpha $.

Recall that E.B.Dynkin classified maximal subgroups of classical groups in
\cite{D}. Since $B_2$ does not occur in his Table 1, 
every irreducible representation of $B_2$ defines a maximal subgroup
by Theorem 1.5 in \cite{D}. In particular,
$H$ is a maximal connected subgroup in $G$.
Therefore, either (i) $H$ is the connected component of $G^\alpha $
or (ii) $\alpha ={\rm Id}$.
Now, (ii) implies that ${\rm Ad}(g_1) = \theta $, thus $\theta $ is an inner automorphism of $G$,
which is not the case. So we are left with (i).
Applying the same argument to $\beta = \alpha ^2$, we see that either (i1) $H$ is the connected component of $G^\beta $
or (i2) $\beta ={\rm Id}$. Since $\beta $ is certainly an inner automorphism, (i1) would imply
that $H$ is the centralizer of an element of $G$. However, all centralizers have even codimension in $G$ and
${\rm codim} (H) = 35$. On the other hand, if (i2) were true then $H$ would be a symmetric subgroup in $G$. The list of
symmetric spaces shows that this
is not the case. The contradiction just obtained completes the proof. 
\end{proof}


\begin{thebibliography}{Name}
\bibitem[A]{A} D.N.Akhiezer, Lie group actions in complex analysis, Vieweg, Braunschweig - Wiesbaden, 1995.
\bibitem[AH]{AH} D. N. Akhiezer, P. Heinzner, Spherical Stein spaces, Manuscripta Math. 114, 
	327-334 (2004)
\bibitem[AV]{AV} D. N. Akhiezer, E. B. Vinberg, Weakly symmetric spaces and 
	spherical varieties,
	Transform. Groups, Vol. 4, No. 1, 1999, 3-24
\bibitem[D]{D} E. B. Dynkin, Maximal subgroups of the classical groups, Amer. Math. Soc. Transl., Ser. 2, vol. 6 (1957), 245 - 378 
\bibitem[FT]{FT} J. Faraut, E. G. F. Thomas, Invariant Hilbert spaces of
	holomorphic functions, J. Lie Theory, Vol. 9 (1999), 383-402
\bibitem[H-C]{H-C} Harish-Chandra, Discrete series for semisimple Lie groups II, Acta Math. 116 
	(1966), 1-111
\bibitem[OV]{OV} A. L. Onishchik, E. B. Vinberg, Lie groups and algebraic groups,
	Springer 1990
\end{thebibliography}
\end{document}